\documentclass{article}

\usepackage{proof}
\usepackage{latexsym}

\usepackage{amssymb}
\usepackage{cite}

\newtheorem{theorem}{Theorem}[section]
\newtheorem{definition}[theorem]{Definition}
\newtheorem{proposition}[theorem]{Proposition}
\newtheorem{lemma}[theorem]{Lemma}
\newtheorem{corollary}[theorem]{Corollary}

\begin{document}

\title{Predicatively computable functions on sets}


\author{
Toshiyasu Arai
\thanks{I'd like to thank Sebastian Eberhard
to point out the
sloppy proofs in section \ref{sec:HF}, and Arnold Beckmann for his interests in this article.}
\\
Graduate School of Science,
Chiba University
\\
1-33, Yayoi-cho, Inage-ku,
Chiba, 263-8522, JAPAN
\\
tosarai@faculty.chiba-u.jp
}

\maketitle

\begin{abstract}
Inspired from a joint work by A. Beckmann, S. Buss and S. Friedman,
we propose a class of set-theoretic functions, predicatively computable set functions.
Each function in this class is polynomial time computable when we restrict to finite binary strings.
\end{abstract}

\section{Introduction}

Bellantoni and Cook\cite{bellantonicook} introduced a class $B$ of functions on finite binary strings.
Arguments of each function $f$ in the class $B$ are divided into \textit{normal} arguments $\vec{x}$ and \textit{safe} arguments $\vec{a}$,
and denoted $f(\vec{x}/\vec{a})$\footnote{Here we follow the notation in \cite{CRM} using slash (/) instead of semicolon (;) to distinguish arguments.} .
Let $\epsilon$ denote the empty string, and $si$ the concatenated string obtained from the binary string $s$ and $i=0,1$.
The class $B$ is generated from initial functions (projections, zero, binary successors $s_{i}(-/s)=si\,(i=0,1)$, 
 the predecessor $p(-/\epsilon)=\epsilon$, $p(-/si)=s$, the conditional(parity test) $C(-/a,b,c)=b$ if $a=s1$, $=c$ otherwise)
 by operating safe composition $f(\vec{x}/\vec{a})=h(\vec{r}(\vec{x}/- )/\vec{t}(\vec{x}/\vec{a}))$
 and predicative recursion on notation $f(\epsilon,\vec{x}/\vec{a})=g(\vec{x}/\vec{a})$ and $f(si,\vec{x}/\vec{a})=h_{i}(s,\vec{x}/\vec{a},f(s,\vec{x}/\vec{a}))$ for $i=0,1$.
It is shown in \cite{bellantonicook} that the polynomial time computable functions are exactly those functions in $B$ having no safe arguments.

It seems to me that the class $B$ not only characterize the class of the polynomial time computable functions,
but also is of foundational importance since each function in $B$ is computable \textit{predicatively}.
By computability we mean that each object reaches to a canonical form by some computations.
However a general concept `computability' involves possibly infinite searches or 
at least the notion of finite computations in general as completed processes.
This is not justified predicatively.
For example a substitution of $f(s,\vec{x}/\vec{a})$ in a normal argument, 
$f(si,\vec{x}/\vec{a})=h_{i}(s,\vec{x},f(s,\vec{x}/\vec{a})/\vec{a})$ is hard to justify predicatively
since it assumes a hypothetical computation of $f(s,\vec{x}/\vec{a})$ to be completed.
On the other side, we see that a computation process of each function $f(\vec{x}/\vec{a})$ in $B$ can be obtained by 
imitating the generating process of normal arguments $\vec{x}$.
In the computation process the safe arguments $\vec{a}$ act only as \textit{names}.
In other words we don't need to know the values (canonical forms) of $\vec{a}$,
but need the values of normal arguments $\vec{x}$ from which we know how the arguments are generated from 
$\epsilon$ by rules $s\mapsto si$.
In this sense the predicative recursion on notation is justifiable predicatively.
This observation was implicit in our joint work\cite{AraiMoser} with G. Moser 
to design a path order POP for computations in $B$.

We now ask how to define predicatively justifiable computations on \textit{sets}?
Contrary to binary strings, there seem no canonical forms of sets even for hereditarily finite sets
unless we assume, e.g., the axiom of constructibility.
Let us approach modestly.
First pick some functions on safe arguments to generate sets such as pairing and unions.
Then applying safe composition and a safe set recursion $f(x,\vec{y}/\vec{a})=h(x,\vec{y}/\vec{a},\{f(z,\vec{y}/\vec{a}) : z\in x\})$
to get a class of functions on sets.
Each set is inductively generated, i.e.,  the epsilon relation $z\in x$ is well founded.
Safe set recursion is close to the idea of predicatively computable functions
since we don't need to know the values of intermediate terms $f(z,\vec{y}/\vec{a}) \, (z\in x)$
to continue the computations of $f(x,\vec{y}/\vec{a})$.
Thus a class {\sf PCSF} of predicatively computable set functions is obtained in section \ref{sec:PCSF}.
The class  {\sf PCSF} is a subclass of the class {\sf SRSF} of safe recursive set functions due to A. Beckmann, S. Buss and S. Friedman \cite{CRM}.
Their joint work motivates ours, and is reported in section \ref{sec:BBF}.

In section \ref{sec:string} it is shown that 
each polynomial time computable function on finite binary strings is in the class {\sf PCSF}, cf. Lemma \ref{lem:ptimePCSF}.
In section \ref{sec:HF}
the size of {\sf PCSF} function $f(\vec{x}/\vec{a})$ is seen to be bounded by a polynomial in the sizes of
normal arguments $\vec{x}$, and 
to depend linearly on the safe arguments $\vec{a}$, cf. Theorem \ref{th:size}.
From this we see readily that each {\sf PCSF} function $f(\vec{x}/-)$ on finite binary strings
 is polynomial time computable, cf. Corollary \ref{th:main}.

\section{Safe recursive set functions}\label{sec:BBF}

A. Beckmann, S. Buss and S. Friedman \cite{CRM} introduced a class
{\sf SRSF} of \textit{safe recursive set functions}.
The class {\sf SRSF}  is obtained from Gandy-Jensen rudimentary set functions on \textit{safe arguments}
by safe composition scheme and predicative set (primitive) recursion scheme a l\`a Bellantoni-Cook.

\begin{description}
\item[(Projection)]
$$
\mbox{proj}^{n,m}_{j}(x_{1},\ldots,x_{n}/ x_{n+1},\ldots, x_{n+m})=x_{j}\, (1\leq j\leq n+m)
.$$

\item[(Difference)]
$$
\mbox{diff}(-/ a,b)=a\setminus b
.$$

\item[(Pair)]
$$
\mbox{pair}(- / a,b)=\{a,b\}
.$$

\item[(Bounded Union)]
$$
f(\vec{x}/ \vec{a},b)=\bigcup_{c\in b}g(\vec{x}/\vec{a},c)
.$$

\item[(Safe Composition)]
$$
f(\vec{x}/\vec{a})=h(\vec{r}(\vec{x}/- )/\vec{t}(\vec{x}/\vec{a}))
.$$

\item[(Predicative Set Recursion)]
$$
f(x,\vec{y}/\vec{a})=h(x,\vec{y}/\vec{a},\{f(z,\vec{y}/\vec{a}) : z\in x\})
.$$
\end{description}

They investigate definability and complexity of safe recursive functions.

\begin{enumerate}
\item
For each $f\in{\sf SRSF}$ there exists a polynomial function $q_{f}$ on ordinals such that
$\mbox{{\rm rank}}(f(\vec{x}/\vec{a}))\leq\max(\mbox{{\rm rank}}(\vec{a}))+q_{f}(\mbox{{\rm rank}}(\vec{x}))$.

\item
A set-theoretic function $f(\vec{x}/-)$ on infinite ranks $\vec{x}$ is in {\sf SRSF}
iff it is $\Sigma_{1}$-definable on $\mbox{{\rm SR}}_{n}(\vec{x}):=L_{\mbox{\footnotesize{\rm rank}}(\vec{x})^{n}}^{\mbox{\footnotesize{\rm TC}}(\vec{x})}$ for an $n<\omega$,
where for ordinals $\alpha$ and sets $x$
$L_{\alpha}^{x}$ denotes the $L$-hierarchy relativized to $x$, and ${\rm TC}(x)$ the transitive closure of $x$.

\item
For each $f\in{\sf SRSF}$ there exists a polynomial function $p_{f}$ such that
\\
$card(\mbox{{\rm TC}}(f(\vec{x}/\vec{a})))\leq card(\mbox{{\rm TC}}(\{\vec{x},\vec{a}\}))^{2^{p(\mbox{\footnotesize{\rm rank}}(\vec{x}))}}$, where $card(x)$ denotes the cardinality of sets $x$.

\item
Under a natural encoding of finite binary strings,
$f\in{\sf SRSF}$ on finite strings are exactly the functions computed by alternating Turing machines running in exponential time with polynomially many alternations.
\end{enumerate}

It seems to me that it is hard to justify the class {\sf SRSF}  predicatively.
The problem lies in {\bf (Bounded Union)} since it requires us to know
\textit{all} of the elements $c$ in the set $b$ in safe argument.
However we don't know its \textit{value}, but only know its \textit{name} of $b$.
Therefore collecting all the elements of sets in safe argument might not be in the idea of predicatively
justifiable computations.

\section{Predicatively computable set functions}\label{sec:PCSF}

Let me propose a subclass {\sf PCSF} of \textit{predicatively} computable set functions.
First a subclass ${\sf PCSF}^{-}$ of {\sf PCSF} is introduced.

Each function $f$ in the subclass ${\sf PCSF}^{-}$ has no normal arguments $f(-/\vec{a})$.
Initial functions in ${\sf PCSF}^{-}$ are {\bf (Projection)} on safe arguments, $\mbox{proj}^{-,m}_{j}(-/a_{1},\ldots,a_{m})=a_{j}$,  {\bf (Pair)},
{\bf (Null)},  
{\bf (Union)},
and
{\bf (Conditional$\in$)}.

\begin{description}

\item[(Null)]
$$
\mbox{null}(-/-)=0=\emptyset
.$$

\item[(Union)]
$$
\mbox{union}(-/a)=\cup a
.$$

\item[(Conditional$\in$)]
$$
\mbox{Cond}_{\in}(-/ a,b,c,d)=\left\{
\begin{array}{ll}
a & \mbox{{\rm if} } c\in d\\
b & \mbox{{\rm otherwise}}
\end{array}
\right.
$$
\end{description}

The class ${\sf PCSF}^{-}$ is closed under composition $f(-/\vec{a})=h(-/\vec{t}(-/\vec{a}))$, and
{\bf (Safe Separation)}.
\begin{description}
\item[(Safe Separation)]
$$
f(-/\vec{a},c)=c\cap\{b: h(-/\vec{a},b)\neq 0\}=\{b\in c: h(-/\vec{a},b)\neq 0\}
.$$
\end{description}
The class {\sf PCSF} is then obtained from ${\sf PCSF}^{-}$ and {\bf (Projection)} $\mbox{proj}^{n,m}_{j}$
by operating {\bf (Safe Composition)} and {\bf (Predicative Set Recursion)}.

A relation $R(\vec{x}/\vec{a})$ is in {\sf PCSF} if its characteristic function $\chi_{R}(\vec{x}/\vec{a})$ is in the class.
($\chi_{R}(\vec{x}/\vec{a})=1$ if 
 $R(\vec{x}/\vec{a})$, $\chi_{R}(\vec{x}/\vec{a})=0$ otherwise.)
 \\
 
 \noindent
 {\bf Remark}.
 It is open, but unlikely the case that the class {\sf PCSF} is closed under the following safe separation scheme.
 $$
f(\vec{x}/\vec{a},c)=c\cap\{b: h(\vec{x}/\vec{a},b)\neq 0\}=\{b\in c: h(\vec{x}/\vec{a},b)\neq 0\}
.$$
 
Recall that a function $f$ is said to be  \textit{simple} iff $R(f(-/\vec{a}),\vec{b})$ is $\Delta_{0}$ for any $\Delta_{0}$-relations $R$.
As in \cite{Jensen} we see the following proposition.

\begin{proposition}\label{prp:simple}
Each $f\in{\sf PCSF}^{-}$ is a simple function. 
Hence $f$ is a $\Delta_{0}$-function in the sense that its graph is $\Delta_{0}$.
\end{proposition}
 
As in \cite{jensenkarp, rathjen} we see the following proposition.
Proposition \ref{prp:rathjen22}.\ref{prp:rathjen22.2a} tells us that a relation is in ${\sf PCSF}^{-}$ 
iff it is rudimentary, cf. \cite{Jensen}.

As in set-theoretic literature,
$b'c=\bigcup\{d: \langle c,d\rangle\in b\}$, which is the unique element $d$ such that $\langle c,d\rangle\in b$
if such a $d$ exists,
and $b''a=\{b'c: c\in a\}$.

\begin{proposition}\label{prp:rathjen22}
\begin{enumerate}
\item\label{prp:rathjen22.-1}
${\rm diff}(-/a,b)=a\setminus b$ is in ${\sf PCSF}^{-}$.

\item\label{prp:rathjen22.0}
If $g(\vec{x}/\vec{a},\vec{b})$ is in {\sf PCSF}, then so is $f$, where
$f(\vec{x},\vec{y}/\vec{b})=g(\vec{x}/\vec{y},\vec{b})$.

\item\label{prp:rathjen22.1}
If $g,h,R$ are in {\sf PCSF}, then so is $f$, where
$f(\vec{x}/\vec{a})=g(\vec{x}/\vec{a})$ if $R(\vec{x}/\vec{a})$, and $f(\vec{x}/\vec{a})=h(\vec{x}/\vec{a})$ else.

\item\label{prp:rathjen22.2}
The class of relations in  {\sf PCSF} is closed under Boolean operations.

\item\label{prp:rathjen22.2a}

A relation $R(-/\vec{a})$ is $\Delta_{0}$ iff its characteristic function $\chi_{R}$ is in ${\sf PCSF}^{-}$.

\item\label{prp:rathjen22.2b}
$f(-/b,c)=b'c=\bigcup\{d\in\cup\cup b:
\langle c, d\rangle\in b\}$ is in ${\sf PCSF}^{-}$
for the $\Delta_{0}$-relation $\langle c, d\rangle\in b$
where
$\langle c,a\rangle:=\{\{c\},\{c,a\}\}$.

\item\label{prp:rathjen22.4}
If $h$ is in {\sf PCSF}, then so is
$f(x,\vec{y}/\vec{a})=h(x,\vec{y}/\vec{a},\bigcup\{f(z,\vec{y}/\vec{a}) : z\in x\})$.

\item\label{prp:rathjen22.3}{\rm (Cf. {\bf (Bounded Union)}.)}

If $h$ is in {\sf PCSF}, then so is $f$, where
$f(x,\vec{y}/\vec{a})=\bigcup\{h(z,\vec{y}/\vec{a}) : z\in x\}$.

\item\label{prp:rathjen22.5}
If $h,R$ are in {\sf PCSF}, then so are $f,g$, where
$f(x,\vec{y}/\vec{a})=\bigcup\{h(z,\vec{y}/\vec{a}) : z\in x,\, R(z,\vec{y}/\vec{a})\}$
and
$g(x,\vec{y}/\vec{a})=\{h(z,\vec{y}/\vec{a}) : z\in x,\, R(z,\vec{y}/\vec{a})\}$.

\item\label{prp:rathjen22.2c}
$\!\upharpoonright\!(x/a)=a\!\upharpoonright\! x=\{\langle z, a'z\rangle:z\in x\}$ and
$rng(x/a)=a''x$ are in ${\sf PCSF}$.

\item\label{prp:rathjen22.7}
The transitive closure ${\rm TC}(x/-)=x\cup\bigcup\{{\rm TC}(y/-):y\in x\}$ and the rank 
${\rm rank}(x/-)=\bigcup\{{\rm rank}(y/-)+1: y\in x\}$ are in {\sf PCSF}.

\item\label{prp:rathjen22.4a}
If $h$ is in {\sf PCSF}, then so is
\begin{description}
\item[(Predicative Function Recursion)]
\[
f(x,\vec{y}/\vec{a})=h(x,\vec{y}/\vec{a},f\!\upharpoonright\! x)
\]
\end{description}
where
$f\!\upharpoonright\! x:=\{\langle z,f(z,\vec{y}/\vec{a})\rangle : z\in x\}$.

Conversely any {\sf PCSF}-function is generated from ${\sf PCSF}^{-}$-functions and {\bf (Projection)} by
{\bf (Safe Composition)} and {\bf (Predicative Function Recursion)}.

\item\label{prp:rathjen22.10}
Let $R$ be a $\Delta_{0}$-relation.
Assume that $\forall x\exists ! y[y\in z \land R(x,y,z,\vec{a})]$.
Let $f(x,z/\vec{a})=y$ iff $y\in z\land R(x,y,z,\vec{a})$.
Then $f$ is in {\sf PCSF}.

\end{enumerate}

\end{proposition}
\textit{Proof}.\hspace{2mm}
\ref{prp:rathjen22}.\ref{prp:rathjen22.-1}.
$a\setminus b=\{c\in a:c\not\in b\}=\{c\in a:{\rm Cond}_{\in}(-/0,1,c,b)\}$ by {\bf (Safe Separation)}.
\\

\noindent
\ref{prp:rathjen22}.\ref{prp:rathjen22.2a}.
If $\chi_{R}\in{\sf PCSF}^{-}$, then $R(\vec{a})\leftrightarrow \chi_{R}(-/\vec{a})=1$ is a $\Delta_{0}$-relation
by Proposition \ref{prp:simple}.

Conversely consider a relation $R(-/\vec{a},c)\equiv \exists b\in c\, Q(-/\vec{a},b)$ with
a $\Delta_{0}$-relation $Q$.
Then $f(-/\vec{a},c)=c\cap \{b: Q(-/\vec{a},b)\}=c\cap\{b:\chi_{Q}(-/\vec{a},b)\neq 0\}$
is in ${\sf PCSF}^{-}$.
Hence so is $\chi_{R}(-/\vec{a},c)=\mbox{Cond}_{\in}(-/ 0,1,0,\{f(-/\vec{a},c)\})$.
For disjunctions $R(-/\vec{a})\lor Q(-/\vec{a})$ use the finite union 
$\chi_{R}(-/\vec{a})\cup\chi_{Q}(-/\vec{a})$, and for negations $R(-/\vec{a})$ use the conditional
${\rm Cond}_{\in}(-/0,1,0,\chi_{R}(-/\vec{a}))$.
\\

\noindent
\ref{prp:rathjen22}.\ref{prp:rathjen22.3}.
Let $g(z,x,\vec{y}/\vec{a},b)=h(z,\vec{y}/\vec{a})$ if $z\in x$, and
$g(z,x,\vec{y}/\vec{a},b)=b$ otherwise, where $z\in x$ is in {\sf PCSF} by 
{\bf (Conditional$\in$)} and Proposition \ref{prp:rathjen22}.\ref{prp:rathjen22.0}.
Let 
\\
$G(u,x,\vec{y}/\vec{a})=g(u,x,\vec{y}/\vec{a},\bigcup\{G(z,x,\vec{y}/\vec{a}) : z\in u\})$.
Then $G$ is in {\sf PCSF} by Proposition \ref{prp:rathjen22}.\ref{prp:rathjen22.4}, and
\begin{eqnarray*}
G(x,x,\vec{y}/\vec{a}) & = & g(x,x,\vec{y}/\vec{a},\bigcup\{G(z,x,\vec{y}/\vec{a}) : z\in x\})
= \bigcup\{G(z,x,\vec{y}/\vec{a}) : z\in x\}
 \\
  & = & \bigcup\{g(z,x,\vec{y}/\vec{a},\bigcup\{G(u,x,\vec{y}/\vec{a}) : u\in z\}) : z\in x\}
=\bigcup\{h(z,\vec{y}/\vec{a}) : z\in x\}
\end{eqnarray*}
\ref{prp:rathjen22}.\ref{prp:rathjen22.5}.
By Proposition \ref{prp:rathjen22}.\ref{prp:rathjen22.3}
$f(x,\vec{y}/\vec{a})=\bigcup\{\mbox{Cond}_{\in}(-/h(z,\vec{y}/\vec{a}),0,0,\chi_{R}(z,\vec{y}/\vec{a})) : z\in x\}$ is in 
{\sf PCSF}.
Then so is $g(x,\vec{y}/\vec{a})=\bigcup\{\{h(z,\vec{y}/\vec{a})\} : z\in x,\, R(z,\vec{y}/\vec{a})\}$.
\\

\noindent
\ref{prp:rathjen22}.\ref{prp:rathjen22.2c}.
By Propositions \ref{prp:rathjen22}.\ref{prp:rathjen22.2b} and \ref{prp:rathjen22}.\ref{prp:rathjen22.5}
both $\!\upharpoonright\!(x/a)=a\!\upharpoonright\! x=\{\langle z, a'z\rangle:z\in x\}$ and $a''x=\bigcup\{a'y : y\in x\}$ are in ${\sf PCSF}$.
\\

\noindent
\ref{prp:rathjen22}.\ref{prp:rathjen22.7}.
Let
$f(x/-)=(\bigcup\{f(y/-): y\in x\})+1$ for $a+1=a\cup\{a\}$.
Then $f(x/-)={\rm rank}(x/-)+1$ and ${\rm rank}(x/-)=\bigcup f(x/-)=\bigcup\{u: u\in f(x/-)\}$
since $a={\rm rank}(x/-)$ is transitive, i.e., $\bigcup a\subset a$.
\\

\noindent
\ref{prp:rathjen22}.\ref{prp:rathjen22.4a}.
Let $k(x,\vec{y}/\vec{a})=u\cup \{\langle z,h(z,\vec{y}/\vec{a},u\!\upharpoonright\! z)\rangle : z\in x\}$
where $u=\bigcup\{k(z,\vec{y}/\vec{a}) : z\in x\}$.
Then $k$ is in {\sf PCSF} by Propostions \ref{prp:rathjen22}.\ref{prp:rathjen22.4}, \ref{prp:rathjen22}.\ref{prp:rathjen22.5}
and \ref{prp:rathjen22}.\ref{prp:rathjen22.2c}.

Suppose
\begin{equation}\label{eq:rathjen22.4a}
k(x,\vec{y}/\vec{a})=\{\langle z,f(z,\vec{y}/\vec{a})\rangle: z\in\mbox{TC}(x/-)\}=f\!\upharpoonright\! \mbox{TC}(x/-)
\end{equation}
Then we have for
$k(x)\!\upharpoonright\! x=\{\langle z, k(x,\vec{y}/\vec{a})'z \rangle :z\in x\}$ and $z\in x$,
$(k(x)\!\upharpoonright\! x)(z)=f(z,\vec{y}/\vec{a})$.
Hence
$f(x,\vec{y}/\vec{a})=h(x,\vec{y}/\vec{a},f\!\upharpoonright\! x)=h(x,\vec{y}/\vec{a},k(x)\!\upharpoonright\! x)$
is in {\sf PCSF}.

It remains to show (\ref{eq:rathjen22.4a}) by induction on $x$.
By IH(=Induction Hypothesis) we have
$k(x,\vec{y}/\vec{a})=\bigcup\{f\!\upharpoonright\! \mbox{TC}(z/-) : z\in x\}\cup\{\langle z,h(z,\vec{y}/\vec{a},f\!\upharpoonright\! z)\rangle : z\in x\}$.
Hence by the definition of $f$ we have
$k(x,\vec{y}/\vec{a})=\bigcup\{f\!\upharpoonright\! \mbox{TC}(z/-) : z\in x\}\cup\{\langle z,f(z,\vec{y}/\vec{a})\rangle : z\in x\}$.
This shows (\ref{eq:rathjen22.4a}), and {\sf PCSF} is closed under {\bf (Predicative Function Recursion)}.

Conversely let $f$ be defined from $h$ by {\bf (Predicative Set Recursion)} as
$f(x,\vec{y}/\vec{a})=h(x,\vec{y}/\vec{a},\{f(z,\vec{y}/\vec{a}) : z\in x\})$.
Then $f(x,\vec{y}/\vec{a})=h(x,\vec{y}/\vec{a},(f\!\upharpoonright\! x)''x)$.
$h_{0}(x,\vec{y}/\vec{a},b)=h(x,\vec{y}/\vec{a}, b''x)$ is in {\sf PCSF} by Proposition \ref{prp:rathjen22}.\ref{prp:rathjen22.2c}.
Hence
$f$ is defined from $h_{0}$ by {\bf (Predicative Function Recursion)}.
\\

\noindent
\ref{prp:rathjen22}.\ref{prp:rathjen22.10}.
By Propositions \ref{prp:rathjen22}.\ref{prp:rathjen22.2} and \ref{prp:rathjen22}.\ref{prp:rathjen22.2a},
$\Delta_{0}$-relation $R(x,y,z,\vec{a})$ defines a relation $R(x,y,z/\vec{a})$ in {\sf PCSF}.
So is $f(x,z/\vec{a})=\bigcup\{y : y\in z,  R(x,y,z/\vec{a})\}$ by Proposition \ref{prp:rathjen22}.\ref{prp:rathjen22.5}.

\hspace*{\fill} $\Box$

\section{Polytime function on finite strings}\label{sec:string}

Let $\mathbb{HF}$ denote the set of all \textit{hereditarily finite sets}.
Let us encode finite (binary) strings by hereditarily finite sets, $\nu:{}^{<\omega}2\to\mathbb{HF}$ 
slightly modified from \cite{CRM}.

$\nu(\epsilon)=0=\emptyset$ ($\epsilon$ is the empty string.)
$\nu(s i)=\langle i+1,\nu(s)\rangle=\{\{i+1\},\{i+1,\nu(s)\}\}\,(i=0,1)$. $1=\{0\}, 2=\{0,1\}$.
For example, $\nu(100)=\langle 1,\langle 1,\langle 2,0\rangle\rangle\rangle$.

\begin{lemma}\label{lem:ptimePCSF}
For each polynomial time computable function $f(\vec{s})$ there exists a function $F$ in {\sf PCSF}
such that for any finite strings $\vec{s}$
\[
F(\nu(\vec{s})/-)=\nu(f(\vec{s}))
.\]
\end{lemma}
\textit{Proof}.\hspace{2mm}
Let $B$ denote the class of safe recursive functions on binary finite  strings in \cite{bellantonicook}.
We show inductively that for each $f(\vec{s}/\vec{a})\in B$ there exists a function $F$ in {\sf PCSF}
such that for any finite strings $\vec{s},\vec{t}$,
$
F(\nu(\vec{s})/\nu(\vec{t}))=\nu(f(\vec{s}/\vec{t}))$.

For the binary successor $s_{i}(-/s)=si\,(i=0,1)$, $S_{i}(-/a)=\{\{i+1\},\{i+1,a\}\}$ does the job.

For the predecessor $p(-/\epsilon)=\epsilon$, $p(-/si)=s$,
first let by Propositions  \ref{prp:rathjen22}.\ref{prp:rathjen22.1} and \ref{prp:rathjen22}.\ref{prp:rathjen22.2a}
$$
\mbox{pred}(-/a)=\left\{
\begin{array}{ll}
\cup a
& \mbox{{\rm if} } \exists b,c\in a[a=\{b,c\}] \\
0 & \mbox{{\rm otherwise}}
\end{array}
\right.
$$
Namely $\mbox{pred}(-/\{b,c\})=b\cup c$.
Then let $P(-/a)=(\mbox{pred}(-/\mbox{pred}(-/a)))\setminus\{0,1\}$ by Proposition \ref{prp:rathjen22}.\ref{prp:rathjen22.-1}.
We have $P(-/\nu(si))=((i+1)\cup\nu(s))\setminus\{0,1\}=\nu(s)$ since $\{0,1\}\cap\nu(s)=\emptyset$.

Next consider conditional(parity test) $C(-/a,b,c)=b$ if $a=s1$, $=c$ otherwise.
Since $2\neq\nu(s)$ and $\{2\}\in\nu(a)\Leftrightarrow a=s1$,
$f(-/a,b,c)=\mbox{Cond}_{\in}(-/b,c,\{2\},a)$ enjoys $f(-/\nu(a),\nu(b),\nu(c))=\nu(C(-/a,b,c))$.

The case when $f(\vec{s}/\vec{a})=h(\vec{r}(\vec{s}/-)/\vec{t}(\vec{s}/\vec{a}))$ is defined from $h,\vec{r},\vec{t}$
by predicative composition is seen from IH.

Finally consider
predicative recursion on notation.
$f(\epsilon,\vec{x}/\vec{a})=g(\vec{x}/\vec{a})$ and $f(si,\vec{x}/\vec{a})=h_{i}(s,\vec{x}/\vec{a},f(s,\vec{x}/\vec{a}))$ for $i=0,1$.
Let $G$ and $H_{i}$ be functions in {\sf PCSF} for $g$ and $h_{i}$, resp.
Define $F$ as follows.
Let $i=0,1$, and $y$ be such that $\{\{i+1,y\},\langle i+1,y\rangle\}\cap\{0,1,2,\langle 1,0\rangle\}=\emptyset$.
Also let $z\not\in\{0,1,2,\langle 1,0\rangle, \{i+1,y\},\langle i+1,y\rangle: i=0,1, y\geq 0\}$.
\begin{eqnarray*}
F(0,\vec{x}/\vec{a}) & := & G(\vec{x}/\vec{a})
\\
F(\langle 1,0\rangle,\vec{x}/\vec{a}) & := & H_{0}(0,\vec{x}/\vec{a},G(\vec{x}/\vec{a}))
\\
F(i+1,\vec{x}/\vec{a}) & := & F(\{i+1\},\vec{x}/\vec{a}):=0
\\
F(\{i+1,y\},\vec{x}/\vec{a})
& := & H_{i}(y,\vec{x}/\vec{a},\bigcup\{F(z,\vec{x}/\vec{a}) : z\in\{i+1,y\}\})
=H_{i}(y,\vec{x}/\vec{a},F(y,\vec{x}/\vec{a}))
\\\
F(\langle i+1,y\rangle,\vec{x}/\vec{a}) & := & \bigcup\{F(z,\vec{x}/\vec{a}) : z\in\langle i+1,y\rangle\}=F(\{i+1,y\},\vec{x}/\vec{a})
\\
F(z,\vec{x}/\vec{a}) & := & 0
\end{eqnarray*}
Then 
$F(\nu(0),\vec{x}/\vec{a})=F(\langle 1,0\rangle,\vec{x}/\vec{a})=H_{0}(0,\vec{x}/\vec{a},G(\vec{x}/\vec{a}))
=H_{0}(0,\vec{x}/\vec{a},F(0,\vec{x}/\vec{a}))$.
Also we compute for $i=0,1$, if $\lnot(s=\epsilon \land i=0)$,
$F(\nu(si),\vec{x}/\vec{a})=F(\langle i+1,\nu(s)\rangle,\vec{x}/\vec{a})=F(\{i+1,\nu(s)\},\vec{x}/\vec{a})
=H_{i}(\nu(s),\vec{x}/\vec{a},F(\nu(s),\vec{x}/\vec{a}))$.
\hspace*{\fill} $\Box$
\\

\noindent
{\bf Remark}.
Lemma \ref{lem:ptimePCSF} holds also for
a subclass $\mbox{{\sf PCSF}}^{\prime}$.
The initial functions in the subclass are
projections $\mbox{proj}^{n,m}_{j}$, $\mbox{diff}(-/ a,b)$,
$S(-/a)=\{a\}$,
$\mbox{pred}(-/a)$ in the proof of Lemma \ref{lem:ptimePCSF}, $\mbox{Cond}_{\in}(-/a,b,c,d)$ 
and $\mbox{finunion}(-/ a,b)=a\cup b$.
The class $\mbox{{\sf PCSF}}^{\prime}$ is closed under {\bf (Safe Composition)} and 
the scheme
$f(x,\vec{y}/\vec{a})=h(x,\vec{y}/\vec{a},\bigcup\{f(z,\vec{y}/\vec{a}) : z\in x\})$, cf. Proposition \ref{prp:rathjen22}.\ref{prp:rathjen22.4}.

Moreover {\bf (Safe Separation)} is needed only in defining
$\mbox{diff}$,
$b'c$ (Proposition \ref{prp:rathjen22}.\ref{prp:rathjen22.2b}) and $\mbox{pred}(-/a)$ for Lemma \ref{lem:ptimePCSF}.
Namely the separation $\mbox{diff}(-/a,b)=\{c\in a:c\not\in b\}$, 
$f(-/b,c,a)=\{d\in a: \langle c,d\rangle\in b\}$, $g(-/a)=\{b\in a: \exists c\in a[a=\{b,c\}]\}$
and $h(-/b,a)=\{c\in a: a=\{b,c\}\}$.

\section{Predicatively computable functions on $\mathbb{HF}$}\label{sec:HF}

Let us restrict our attention to hereditarily finite sets $\mathbb{HF}$.
$X,Y,Z,U,\ldots,A,B$ denote hereditarily finite sets.
Each function $f$ in {\sf PCSF} is a function on $\mathbb{HF}$ when it is restricted to $\mathbb{HF}$.

The size of $f(\vec{x}/\vec{a})$ is seen to be bounded by a polynomial in the sizes of
normal arguments $\vec{x}$, and 
depend \textit{linearly} on the safe arguments $\vec{a}$, cf. Theorem \ref{th:size}.
This readily yields  the converse of Lemma \ref{lem:ptimePCSF}, cf. Corollary \ref{th:main}.

For a polynomial $p(\vec{x})$ and hereditarily finite sets 
$\vec{X}=X_{1},\ldots,X_{n}\in\mathbb{HF}$, put
\begin{eqnarray*}
cT(X) & := & card(\mbox{{\rm TC}}(X))
\\
pt(\vec{X}) & := & p(cT(X_{1}),\ldots,cT(X_{n}))
\end{eqnarray*}

A polynomial $p(x_{1},\ldots,x_{n})$ is said to be \textit{weakly monotonic}
if
 $\forall i\leq n(x_{i}\leq y_{i})\Rightarrow p(x_{1},\ldots,x_{n})\leq p(y_{1},\ldots,y_{n})$.
 
\begin{theorem}\label{th:size}
For each (definition of) function $f(x_{1},\ldots,x_{n}/a_{1},\ldots,a_{m})\in{\sf PCSF}$
there exists a weakly monotonic polynomial
$p_{f}(\vec{x})$ such that for any hereditarily finite sets $\vec{X}=X_{1},\ldots,X_{n}$ and $\vec{A}=A_{1},\ldots,A_{m}$,
the size of the set difference of the transitive closures of $f(\vec{X}/\vec{A})$ and of $\cup S(\vec{A})$
is bounded by $pt_{f}(\vec{X})$:
\[
card(
\mbox{{\rm TC}}(
f(\vec{X}/\vec{A}))
\setminus
\mbox{{\rm TC}}(
\cup S(\vec{A}))
)
\leq
pt_{f}(\vec{X})
\]
and
\[
cT(f(\vec{X}/\vec{A}))\leq pt_{f}(\vec{X})+cT(\cup S(\vec{A}))
\]
where
$pt_{f}(\vec{X}):=p_{f}(cT(X_{1}),\ldots,cT(X_{n}))$ and $\cup S(\vec{A}) :=S(A_{1})\cup\cdots\cup S(A_{m})$
with $S(A)=A\cup\{A\}$.
\end{theorem}

The theorem says that safe arguments $\vec{A}$ are never duplicated.

\begin{corollary}\label{cor:times}
The Cartesian product $\mbox{{\rm prod}}(-/a,b)=a\times b$ is not in {\sf PCSF}.
Even $f(-/a)=\{0\}\times a=\{\langle 0,b\rangle: b\in a\}\not\in{\sf PCSF}$.

On the other side, $f(x,y/-)=x\times y$ is in {\sf PCSF}.
\end{corollary}
\textit{Proof}.\hspace{2mm}
Consider the hereditarily finite  sets $a_{n}=\{2,\ldots,n\}$ for $n\geq 2$.
Then $\langle 0,b\rangle, \{0,b\}\not\in {\rm TC}(a_{n})$ for any $b\in a_{n}$, and 
$cT(\{0\}\times a_{n})\geq  cT(a_{n})+card(a_{n})$.

On the other hand we have
$x\times y=\bigcup_{u\in x}\bigcup_{v\in y}\{\langle u,v\rangle\}$.
\hspace*{\fill} $\Box$
\\

\noindent
Let us introduce some abbreviations to state and shorten the proof of the following lemma.
For hereditarily finite sets $\{X_{i}, Z_{i}: 1\leq i\leq n\}\cup\{A_{i}:1\leq i\leq k\}\subset \mathbb{HF}$,
let us denote
$\vec{X}=X_{1},\ldots,X_{n}$, $\vec{Z}=Z_{1},\ldots,Z_{n}$,
$\vec{A}=A_{1},\ldots,A_{k}$, 
$\cup S(\vec{A})=S(A_{1})\cup\cdots\cup S(A_{k})$,
$S(A)=A\cup\{A\}$,
and $\vec{X}\in \mbox{{\rm TC}}(\vec{Z}):\Leftrightarrow \forall i\leq n[X_{i}\in \mbox{{\rm TC}}(Z_{i})]$.
Also let
$
\{\vec{g}(\vec{X}/\vec{A}):
\vec{X}\in\mbox{{\rm TC}}(\vec{Z})\}
:=
\{g_{j}(\vec{X}/\vec{A}): \vec{X}\in\mbox{{\rm TC}}(\vec{Z}), 1\leq j\leq m\}
$ for sequences $\vec{g}=g_{1},\ldots,g_{m}$ of functions.

\begin{lemma}\label{prp:coverp}
For each (definition of) function $f(\vec{x}/\vec{b})\in{\sf PCSF}$ with 
$\vec{x}=x_{1},\ldots,x_{n}, \vec{b}=b_{1},\ldots,b_{m}$
there exists a weakly monotonic polynomial
$q_{f}(\vec{x})$ for which the following hold.

For any list $\vec{Z}=Z_{1},\ldots,Z_{n}$ of  hereditarily finite sets $Z_{i}\in\mathbb{HF}$,
any list of functions $\vec{g}(\vec{x}/\vec{a})=g_{1}(\vec{x}/\vec{a}),\ldots,g_{m}(\vec{x}/\vec{a})$ of $g_{i}\in{\sf PCSF}$,
any list $\vec{A}$ of $A_{i}\in\mathbb{HF}$,
the cardinality of the following set (difference) is at most $qt_{f}(\vec{Z})$:
\[
\mbox{{\rm TC}}(
\{f(\vec{X}/\vec{g}(\vec{X}/\vec{A})):
\vec{X}\in \mbox{{\rm TC}}(\vec{Z})\})
\setminus
\mbox{{\rm TC}}(
\cup S(\vec{A})\cup
\{\vec{g}(\vec{X}/\vec{A}): \vec{X}\in \mbox{{\rm TC}}(\vec{Z})\}
)
.\]
\end{lemma}

Lemma \ref{prp:coverp} yields Theorem \ref{th:size} as follows.
For $f\in{\sf PCSF}$ and $g_{j}(\vec{x}/\vec{a})=a_{j}$, i.e., the projection $g_{j}=\mbox{proj}^{n,m}_{n+j}$, 
we have a polynomial 
$q_{f}$ such that for any lists of hereditarily finite sets $\vec{Z}$, $\vec{A}$,
\[
card(\mbox{{\rm TC}}(\{f(\vec{X}/\vec{A}):\vec{X}\in \mbox{{\rm TC}}(\vec{Z})\})
\setminus \mbox{{\rm TC}}(\cup S(\vec{A})))
\leq
qt_{f}(\vec{Z})
.\]
Let $Z_{i}=\{X_{i}\}$.
Then $\mbox{TC}(Z_{i})=\{X_{i}\}\cup \mbox{TC}(X_{i})$ and
$card(
\mbox{{\rm TC}}(
f(\vec{X}/\vec{A}))
\setminus
\mbox{{\rm TC}}(
\cup S(\vec{A}))
)
\leq
pt_{f}(\vec{X})$
for $p_{f}(x_{1},\ldots,x_{n})=q_{f}(x_{1}+1,\ldots,x_{n}+1)$.
\\

\noindent
\textit{Proof} of Lemma \ref{prp:coverp}.\hspace{2mm}
Let us define a natural number $o(f)<\omega$ for each function $f\in{\sf PCSF}$ as follows.
First $o(f)=0$ if $f$ is one of null, pair, projections $\mbox{proj}^{n,m}_{j}$, union, $\mbox{Cond}_{\in}$ and 
functions defined by {\bf (Safe Separation)}.
Second $o(f)=1+\max\{o(h), o(r_{i}), o(t_{j}): i=1,\ldots,n, j=1,\ldots, m\}$
if $f$ is defined by {\bf (Safe Composition)} from $h,r_{1},\ldots,r_{n}$,
$t_{1}\ldots,t_{m}$.
Third $o(f)=1+o(h)$ if $f$ is defined by {\bf (Predicative Set Recursion)} from $h$.
The lemma is shown by induction on the number $o(f)$ assigned to the definition of $f$.
\\
{\bf (Null)}
If $f$ is $\mbox{null}(-/-)=\emptyset$, then $q_{f}(-)=0$.
\\
{\bf (Projection)}
If $f$ is a projection $\mbox{proj}^{n,m}_{i}$,
then 
$f(\vec{X}/\vec{g}(\vec{X}/\vec{A}))$ is one of $X_{i}$ or $g_{i-n}(\vec{X}/\vec{A})$.
In the former case $q_{f}(\vec{x})=x_{i}$,
while in the latter case $q_{f}(\vec{x})=0$.
\\
{\bf (Pair)}
If $f$ is the pair $\mbox{pair}(-/A_{1},A_{2})=\{A_{1},A_{2}\}$, then $q_{f}(\vec{x})=1$.
\\
{\bf (Union)}
If $f$ is the union $\mbox{union}(-/A_{1})=\cup A_{1}$, then $q_{f}(\vec{x})=1$.
\\
{\bf (Conditional$\in$)}
If $f$ is the conditional $\mbox{Cond}_{\in}(-/A_{1},A_{2},A_{3},A_{4})\in\{A_{1},A_{2}\}$, then
$q_{f}(\vec{x})=0$.
\\
{\bf (Safe Separation)}
If $f$ is defined from $h$ by
{\bf (Safe Separation)} $f(-/\vec{A},C)=\{B\in C: h(-/\vec{A},B)\neq 0\}\subset C$,
then $q_{f}(\vec{x})=1$.
\\
{\bf (Safe Composition)}
\\
Consider the case when $f$ is defined from $h$, $\vec{r}$ and $\vec{t}=t_{1},\ldots,t_{k}$ by {\bf (Safe Composition)},
$f(\vec{X}/\vec{g}(\vec{X}/\vec{A}))=
h(\vec{r}(\vec{X}/-)/\vec{t}(\vec{X}/\vec{g}(\vec{X}/\vec{A})))$, where
each $t_{i}(\vec{X}/\vec{g}(\vec{X}/\vec{A}))$ is a {\sf PCSF}-function.

By IH we have a weakly monotonic polynomial $q_{h}(\vec{u},\vec{x})$ such that
for any $\vec{Z},\vec{A}$ and any $\vec{U}$, 
\[
card(\mbox{{\rm TC}}(\mathcal{S}_{h(\vec{t}(\vec{g}))}(\vec{U},\vec{Z}/\vec{A}))\setminus 
\mbox{{\rm TC}}(\cup S(\vec{A})\cup\mathcal{S}_{\vec{t}(\vec{g})}(\vec{Z}/\vec{A})))
\leq qt_{h}(\vec{U},\vec{Z})
\]
where
\begin{eqnarray*}
\mathcal{S}_{h(\vec{t}(\vec{g}))}(\vec{U},\vec{Z}/\vec{A}) & = &
\{h(\vec{Y}/\vec{t}(\vec{X}/\vec{g}(\vec{X}/\vec{A}))) :
\vec{Y}\in\mbox{{\rm TC}}(\vec{U}), \vec{X}\in\mbox{{\rm TC}}(\vec{Z})\}
\\
\mathcal{S}_{\vec{t}(\vec{g})}(\vec{Z}/\vec{A})) & = &
\{\vec{t}(\vec{X}/\vec{g}(\vec{X}/\vec{A})) :
\vec{X}\in\mbox{{\rm TC}}(\vec{Z})\}
\end{eqnarray*}

On the other hand we have a polynomial $q_{t_{i}}$ for each $i=1,\ldots,k$ such that
the size of the following set is bounded by $qt_{t_{i}}(\vec{Z})$:
\[
\mbox{{\rm TC}}(
\{t_{i}(\vec{X}/\vec{g}(\vec{X}/\vec{A})):\vec{X}\in\mbox{{\rm TC}}(\vec{Z})\}
)
\setminus
\mbox{{\rm TC}}(
\cup S(\vec{A})\cup
\{\vec{g}(\vec{X}/\vec{A}):\vec{X}\in\mbox{{\rm TC}}(\vec{Z})\})
\]
Hence for $q_{\vec{t}}(\vec{x})=\sum_{i}q_{t_{i}}(\vec{x})$, $qt_{\vec{t}}(\vec{Z})$ gives an upper bound of the size of the 
following set:
\[
\mbox{{\rm TC}}(
\{\vec{t}(\vec{X}/\vec{g}(\vec{X}/\vec{A})):\vec{X}\in\mbox{{\rm TC}}(\vec{Z})\}
)
\setminus
\mbox{{\rm TC}}(
\cup S(\vec{A})\cup
\{\vec{g}(\vec{X}/\vec{A}):\vec{X}\in\mbox{{\rm TC}}(\vec{Z})\})
\]
Moreover by IH we have a polynomial $q_{r_{i}}(\vec{x})$ such that
$cT(\{r_{i}(\vec{X}/-): \vec{X}\in\mbox{{\rm TC}}(\vec{Z})\})\leq qt_{r_{i}}(\vec{Z})$.
Let 
$\vec{U}=\{r_{1}(\vec{X}/-): \vec{X}\in\mbox{{\rm TC}}(\vec{Z})\},\ldots,\{r_{v}(\vec{X}/-): \vec{X}\in\mbox{{\rm TC}}(\vec{Z})\}$ 
for $\vec{r}=r_{1},\ldots,r_{v}$.
Then
$
\{f(\vec{X}/\vec{g}(\vec{X}/\vec{A})):\vec{X}\in\mbox{{\rm TC}}(\vec{Z})\}
=
\{h(\vec{r}(\vec{X}/\vec{A})/\vec{t}(\vec{X}/\vec{g}(\vec{X}/\vec{A})))
:\vec{X}\in\mbox{{\rm TC}}(\vec{Z})\}
$ is a subset of
$\{h(\vec{Y}/\vec{t}(\vec{X}/\vec{g}(\vec{X}/\vec{A})))
: \vec{Y}\in\mbox{{\rm TC}}(\vec{U}), \vec{X}\in\mbox{{\rm TC}}(\vec{Z})\}
$.
Therefore
\begin{eqnarray*}
&&
card(
\mbox{{\rm TC}}(
\{ f(\vec{X}/\vec{g}(\vec{X}/\vec{A})): \vec{X}\in\mbox{{\rm TC}}(\vec{Z})\})
\setminus
\mbox{{\rm TC}}(
\cup S(\vec{A})\cup
\{\vec{g}(\vec{X}/\vec{A}):\vec{X}\in\mbox{{\rm TC}}(\vec{Z})\})
)
\\
& \leq &
qt_{h}(\vec{U},\vec{Z})+qt_{\vec{t}}(\vec{Z})
\leq
qt_{f}(\vec{Z})
\end{eqnarray*}
where in
$q_{f}(\vec{x})=q_{h}(q_{r_{1}}(\vec{x}),\ldots,q_{r_{v}}(\vec{x}),\vec{x})+q_{\vec{t}}(\vec{x})$, 
$q_{r_{i}}(\vec{x})$ is substituted for each variable $u_{i}$ in $q_{h}(u_{1},\ldots,u_{v},\vec{x})$.
\\
{\bf (Predicative Set Recursion)}
\\
Consider the case when $f$ is defined from $h$ by {\bf (Predicative Set Recursion)},
$f(y,\vec{x}/\vec{a})=h(y,\vec{x}/\vec{a},\{f(z,\vec{x}/\vec{a}): z\in y\})$.
By Proposition \ref{prp:rathjen22}.\ref{prp:rathjen22.5} there exists a {\sf PCSF}-function
$F(Y,\vec{X}/\vec{A})=\{f(Z,\vec{X}/\vec{g}(Z,\vec{X}/\vec{A})): Z\in Y\}$.
Let
\begin{eqnarray*}
\mathcal{D}_{f}(W,\vec{Z}/\vec{A})
& = &
\mbox{{TC}}(\mathcal{S}_{f}(W,\vec{Z}/\vec{A}))
\setminus
\mbox{{TC}}(\cup S(\vec{A})\cup
\mathcal{S}_{\vec{g}}(W,\vec{Z}/\vec{A})
)
\\
\mathcal{S}_{f}(W,\vec{Z}/\vec{A}) 
& = &
\{f(Y,\vec{X}/\vec{g}(Y,\vec{X}/\vec{A})):
Y\in\mbox{{\rm TC}}(W), \vec{X}\in\mbox{{\rm TC}}(\vec{Z})\}
\\
& = &
\{h(Y,\vec{X}/\vec{g}(Y,\vec{X}/\vec{A}),F(Y,\vec{X}/\vec{A})) :
Y\in\mbox{{\rm TC}}(W), \vec{X}\in\mbox{{\rm TC}}(\vec{Z})\}
\\
\mathcal{S}_{\vec{g}}(W,\vec{Z}/\vec{A}) & = &
\{\vec{g}(Y,\vec{X}/\vec{A}) :
Y\in\mbox{{\rm TC}}(W), \vec{X}\in\mbox{{\rm TC}}(\vec{Z})\}
\\
\mathcal{S}_{F}(W,\vec{X}/\vec{A})
& = &
\{F(Y,\vec{X}/\vec{A}): Y\in\mbox{{\rm TC}}(W), \vec{X}\in\mbox{{\rm TC}}(\vec{Z})\}
\end{eqnarray*}
By IH we have a weakly monotonic polynomial $q_{h}(w,\vec{x})$ such that
for any $\vec{Z},\vec{A}$ and any $W$, 
\[
card(\mbox{{\rm TC}}(\mathcal{S}_{f}(W,\vec{Z}/\vec{A}))\setminus 
\mbox{{\rm TC}}(\cup S(\vec{A})\cup\mathcal{S}_{\vec{g}}(W,\vec{Z}/\vec{A})\cup\mathcal{S}_{F}(W,\vec{X}/\vec{A})
))
\leq
qt_{h}(W,\vec{Z})
\]
We have
$card(\mathcal{S}_{F}(W,\vec{X}/\vec{A}))\leq 
cT(W)\prod cT(\vec{Z})$ for
$\prod cT(\vec{Z})=\prod_{i=1}^{n}cT(Z_{i})$.
Hence
\[
card(
\mbox{{\rm TC}}(
\mathcal{S}_{F}(W,\vec{X}/\vec{A}))
\setminus
\mbox{{TC}}(\cup S(\vec{A})\cup
\mathcal{S}_{\vec{g}}(W,\vec{Z}/\vec{A}))
)
\leq
cT(W)\prod cT(\vec{Z})+
card(
\mathcal{D}_{f}(\cup W,\vec{Z}/\vec{A})
)
\]
and
\[
card(
\mathcal{D}_{f}(W,\vec{Z}/\vec{A})
)
\leq
qt_{h}(W,\vec{Z})+cT(W)\prod cT(\vec{Z})+
card(
\mathcal{D}_{f}(\cup W,\vec{Z}/\vec{A})
)
\]
For $\ell\in\omega$, define $\cup^{(\ell)}W$ recursively by $\cup^{(0)}W=W$ and $\cup^{(\ell+1)}W=\cup(\cup^{(\ell)}W)$.
We see inductively that 
\[
card(
\mathcal{D}_{f}(W,\vec{Z}/\vec{A})
)
\leq
\sum_{i<\ell}qt_{h}(\cup^{(i)}W,\vec{Z})+(\sum_{i<\ell}cT(\cup^{(i)}W))\prod cT(\vec{Z})+
card(
\mathcal{D}_{f}(\cup^{(\ell)} W,\vec{Z}/\vec{A})
)
\]
Then for
$\ell=\mbox{rank}(W)\leq cT(W)$, we have $\mbox{{\rm TC}}(\cup^{(\ell)}W)=\emptyset$, and we obtain
\[
card(
\mathcal{D}_{f}(W,\vec{Z}/\vec{A})
)
\leq
\sum_{i<\ell}qt_{h}(\cup^{(i)}W,\vec{Z})+(\sum_{i<\ell}cT(\cup^{(i)}W))\prod cT(\vec{Z})
\]
Therefore
$q_{\vec{f}}(y,\vec{x})=y\cdot q_{\vec{h}}(y,\vec{x})+y^{2}\prod\vec{x}$ works for $f$.

This completes a proof of Lemma \ref{prp:coverp}, and hence of Theorem \ref{th:size}.
\hspace*{\fill} $\Box$

\subsection{Computing on directed acyclic graphs}\label{subsec:dag}
Now we show that any function $f\in{\sf PCSF}$ is polynomial time computable
when we restrict $f$ to $\mathbb{HF}$.
To be specific, let us encode hereditarily finite sets first by DAG's (Directed Acyclic Graphs),
and then encode DAG's by natural numbers.

\begin{definition}\label{df:DAG}
{\rm A} DAG with root {\rm is a triple} $G=(V,E,r)$ {\rm of non-empty finite set} $V$ {\rm of natural numbers},
$E\subset V\times V$ {\rm and} $r\in V$ {\rm such that}
\begin{enumerate}
\item\label{df:DAG1}
{\rm The only node of indegree zero  is} $r${\rm , i.e.,}
$\lnot\exists a\in V[(a,r)\in E]$ {\rm and} $\forall a\in V\setminus\{r\}\exists b\in V[(b,a)\in E]$.

\item\label{df:DAG2}
$\forall (a,b)\in E[a>b]$.
\end{enumerate}
\end{definition}
In what follows a DAG with root is simply said to be a DAG.
$(a,b)\in E$ designates that there is an edge from $a$ to $b$.
From the condition (\ref{df:DAG2}) in Definition \ref{df:DAG} we see that 
$G$ is acyclic.
For a DAG $G=(V,E,r)$ we write $V=V_{G}$, $E=E_{G}$ and $r=r_{G}$.

For nodes $a\in G$, $G|a$ denotes a DAG $G|a=(V_{G}|a,E_{G}|a,a)$
defined by $E_{G}|a=E_{G}\cap(V_{G}|a\times V_{G}|a)$, and for $b\in V_{G}$,
$b\in V_{G}|a$ iff there exists a path from $a$ to $b$ in $G$, i.e., 
there is a sequence $\{(a_{i},b_{i})\}_{i\leq n}\subset E_{G}$ such that
$a_{0}=a$, $b_{n}=b$ and $\forall i<n(b_{i}=a_{i+1})$.

The \textit{rank} $rk_{G}(a)$ of nodes $a$ in $G$ is defined by
$rk_{G}(a)=\max\{rk_{G}(b)+1: (a,b)\in E_{G}\}$, where
$\max\emptyset:=0$.
Then the rank of $G$ is defined by $rk(G)=rk_{G}(r)$.
While the \textit{length} $\ell_{G}(a)$ of the longest path from $r$ to $a$ is defined by
$\ell_{G}(a)=\max\{\ell_{G}(b)+1: (b,a)\in E_{G}\}$.

Since DAG is similar to term graph, we follow terminology in \cite{tgr}.

\begin{definition}\label{bisimilar}
{\rm Let} $G=(V_{G},E_{G},r_{G}),H=(V_{H},E_{H},r_{H})$ {\rm be DAG's.}
\begin{enumerate}
\item
{\rm Each node}
$a\in G$ {\rm encodes a hereditarily finite set} $set_{G}(a)$ {\rm defined by recursion on ranks} $rk_{G}(a)${\rm :}
\[
set_{G}(a)=\{set_{G}(b): (a,b)\in E_{G}\}
.\]
{\rm DAG} $G$ {\rm encodes a hereditarily finite set} $set(G)=set_{G}(r_{G})$.

\item
$a\in G$ {\rm and} $b\in H$ {\rm are} \textit{bisimilar} {\rm (with respect to} $G,H${\rm ), denoted} 
$a\simeq_{G,H}b$ {\rm or simply} $a\simeq b$
{\rm iff} 
$set_{G}(a)=set_{H}(b)$.

$G$ {\rm and} $H$ {\rm are} \textit{bisimilar}{\rm , denoted} $G\simeq H$ {\rm iff}
$r_{G}\simeq_{G,H}r_{H}${\rm , i.e.,} $set(G)=set(H)$.

\item
$G$ {\rm is} \textit{fully collapsed} {\rm iff for any nodes} $a,b$ {\rm in} $G${\rm , if}
$set_{G}(a)=set_{G}(b)$ {\rm then} $a=b$.
\end{enumerate}
\end{definition}
Clearly if $a\simeq_{G}b$, then $rk(a)_{G}=rk_{G}(b)$.

We assume a feasible encoding of finite sequences of natural numbers.
$\langle a_{0},\ldots,a_{n-1}\rangle$ denotes the code of sequence $(a_{0},\ldots,a_{n-1})$ of natural numbers $a_{i}$.
$\lceil G\rceil\in\omega$ denotes the code of DAG $G=(V,E,r)$.
Specifically $\lceil(V,E,r)\rceil=\langle\lceil V\rceil,\lceil E\rceil,r\rangle$,
where for nodes $V=\{r=a_{0}>a_{1}>\cdots>a_{m-1}\}$, its code
$\lceil V\rceil=\langle a_{0},\ldots,a_{m-1}\rangle$,
and
for edges $E=\{e_{0},\ldots,e_{n-1}\}$,
$\lceil E\rceil=\langle\lceil e_{0}\rceil,\ldots,\lceil e_{n-1}\rceil\rangle$,
where $\lceil (a,b)\rceil=\langle a,b\rangle$ and $\lceil e_{0}\rceil>\cdots>\lceil e_{n-1}\rceil$.

It is plain to see that to be a code of a DAG is polynomial time decidable,
and ranks $rk_{G}(a)$ and lengths $\ell_{G}(a)$
of nodes $a$ in $G$ are polynomial time computable from $n=\lceil G\rceil$ and $a$.
Moreover given a code 
$\lceil G\rceil$ 
of a DAG $G$ and a node $a\in V_{G}$, one can compute the code $\lceil G|a\rceil$
in polynomial time.
Therefore let us identify DAG $G$ with its code $\lceil G\rceil$, and, e.g., say that
$G|a$ is polynomial time computable.

Let $|n|=\lfloor\log_{2}(n+1)\rfloor$.
There is a constant $\alpha$ such that for any DAG $G$
\[
cT(set(G))\leq card(V_{G})-1\leq |\lceil G\rceil|\leq \alpha |r_{G}| \cdot card(V_{G})^{2}
\]
and if $G$ is fully collapsed,
\[
cT(set(G))= card(V_{G})-1\leq |\lceil G\rceil|\leq \alpha |r_{G}| \cdot cT(set(G))^{2}
.\]

We say that $G$ is \textit{balanced} if $a\leq card(V_{G|a})$ for any $a\in V_{G}$.
For balanced and fully collapsed
DAG $G$,
$cT(set(G))$ is polynomially related to $|\lceil G\rceil|$.

\begin{proposition}\label{prp:bisimilarfullcollapsed}
\begin{enumerate}
\item\label{prp:bisimilarfullcollapsed1}
Bisimilarity in DAG's is polynomial time decidable.

\item\label{prp:bisimilarfullcollapsed15}
There is a polynomial time function $R$ such that for any given DAG $G$,
$R(G)$ and $G$ are bisimilar and $R(G)$ is balanced with $\lceil R(G)\rceil\leq\lceil G\rceil$.
Moreover if $G$ is fully collapsed, then so is $R(G)$.

\item\label{prp:bisimilarfullcollapsed2}
There is a polynomial time function $c$ such that for any given 
DAG's $G_{0},\ldots,G_{n-1}$,
$c(G_{0},\ldots,G_{n-1})$ is a fully collapsed DAG such that
\[
set(c(G_{0},\ldots,G_{n-1}))=\{set(G_{i}) : i<n\}
.\]
\end{enumerate}
\end{proposition}
\textit{Proof}.\hspace{2mm}
\ref{prp:bisimilarfullcollapsed}.\ref{prp:bisimilarfullcollapsed1}.
Let $b\in_{G}a$ iff there exists an edge $(a,b)\in E_{G}$.
Then
$G\simeq H$ iff 
$\forall a\in_{G}r_{G}\exists b\in_{H}r_{H}(G|a\simeq H|b) \,\&\, \forall b\in_{H}r_{H}\exists a\in_{G}r_{G}(G|a\simeq H|b)$.
A bisimilarity test is performed at most $card(V_{G})\cdot card(V_{H})$ times.
\\

\noindent
\ref{prp:bisimilarfullcollapsed}.\ref{prp:bisimilarfullcollapsed2}.
We can assume that sets $V_{G_{i}}$ are disjoint, for otherwise 
replace $G_{i}$ by $\{i\}\times G_{i}$, where
$V_{\{i\}\times G_{i}}=\{\pi(i,a):a\in V_{G_{i}}\}$ and
$E_{\{i\}\times G_{i}}=\{(\pi(i,a),\pi(i,b)):(a,b)\in E_{G_{i}}\}$ 
for the bijective pairing $\pi(i,j)=\frac{(i+j)(i+j+1)}{2}+j$.
Note that $a>b\Rightarrow \pi(i,a)>\pi(i,b)$.
Let $r=\max\{r_{G{i}}:i<n\}+1$, and $G$ be the joined DAG.
$V_{G}=\{r\}\cup\bigcup_{i<n}V_{G_{i}}$, $r_{G}=r$ and
$E_{G}=\{(r,r_{G_{i}}) : i<n\}\bigcup_{i<n}E_{G_{i}}$.
Clearly $set(G)=\{set(G_{i}) : i<n\}$.

By recursion on ranks define DAG's $\{H_{i}\}_{-1\leq i\leq rk(G)}$ so that
each $H_{i}\simeq G$ and any bisimilar pair $a\simeq_{H_{i}}b$
has ranks larger than $i$, $rk_{H_{i}}(a)=rk_{H_{i}}(b)>i$, as follows.
Let $H_{-1}=G$.
Assume that $H_{i-1}$ has been defined.
Consider $a\in H_{i-1}$ of rank $i$ and its bisimilar class 
$B_{i}(a)=\{b\in V_{H_{i-1}}: b\simeq_{H_{i-1}}a\}$,
and let us share nodes in $B_{i}(a)$.
Note that for $b,c\in B_{i}(a)$ and any $d$, $(b,d)\in E_{H_{i}}\Leftrightarrow (c,d)\in E_{H_{i}}$
by the construction.
Let $a_{i}=\min B_{i}(a)$.
Delete every nodes in $B_{i}(a)$ except $a_{i}$, and
each edge $(d,b)\in E_{H_{i-1}}$ for $b\in B_{i}(a)$
is switched to a new edge $(d,a_{i})$, where $d>b\geq a_{i}$.
The switchings are performed for each $a\in H_{i-1}$ of rank $i$.
The resulting DAG $H_{i}$ is bisimilar to $H_{i-1}$, 
and $a\simeq_{H_{i}}b \Rightarrow rk_{H_{i}}(a)>i$.

Thus $c(G_{0},\ldots,G_{n-1})=H_{rk(G)}$ is fully collapsed and bisimilar to $G$.
\hspace*{\fill} $\Box$
\\

Each $f\in{\sf PCSF}$ on $\mathbb{HF}$ is a polynomial time computable function in the following sense.

\begin{theorem}\label{lem:hfptime}
For each $f\in{\sf PCSF}$, there is a polynomial time computable function $F$ such that
for any balanced and fully collapsed DAG's $\vec{G}$, $\vec{H}$,
$F(\lceil\vec{G}\rceil,\lceil\vec{H}\rceil)$ is a code $\lceil K\rceil$ of a balanced and fully collapsed DAG $K$ such that
$f(set(\vec{G})/ set(\vec{H}))=set(K)$.
\end{theorem}
\textit{Proof}.\hspace{2mm}
This is seen by construction of $f\in{\sf PCSF}$.
We assume that any DAG is transformed to a balanced one if necessary
by Proposition \ref{prp:bisimilarfullcollapsed}.\ref{prp:bisimilarfullcollapsed15}.

{\bf (Pair)} The case when $f$ is the pairing {\sf pair}
follows from Proposition \ref{prp:bisimilarfullcollapsed}.\ref{prp:bisimilarfullcollapsed2}.

{\bf (Union)} 
For DAG $G$, a DAG $H$ such that $set(H)=\cup(set(G))$ is obtained by $r_{H}=r_{G}$,
$V_{H}=\{a\in V_{G}: \ell_{G}(a)\neq 1\}$ and for $a,b\in V_{H}$,
$(a,b)\in E_{H}$ iff either $(a,b)\in E_{G}$ or there is a $c\in V_{G}$ such that
$\ell_{G}(c)=1$ and $(a,c),(c,b)\in E_{G}$.

{\bf (Conditional$\in$)} follows from Proposition \ref{prp:bisimilarfullcollapsed}.\ref{prp:bisimilarfullcollapsed1},
and {\bf (Safe Separation)} follows from IH.

Next consider {\bf (Safe Composition)}
$$
f(\vec{x}/\vec{a})=h(\vec{r}(\vec{x}/- )/\vec{t}(\vec{x}/\vec{a}))
.$$
If all of $h$, $\vec{r}$ and $\vec{t}$ are polynomial time computable on DAG's,
then so is $f$.

Finally consider {\bf (Predicative Set Recursion)}
$$
f(x,\vec{y}/\vec{a})=h(x,\vec{y}/\vec{a},\{f(z,\vec{y}/\vec{a}) : z\in x\})
.$$
Assume that $x,\vec{y},\vec{a}$ are hereditarily finite sets $set(G), set(\vec{H}), set(\vec{K})$ for fully collapsed DAG's $G,\vec{H},\vec{K}$.
Let us describe informally a polynomial time computation of
a fully collapsed DAG $L$ such that $set(L)=f(set(G),set(\vec{H})/set(\vec{K}))$.
By recursion on ranks $rk_{G}(a)$ of nodes $a$ in `circuit' $G$,
assign a DAG $L_{a}$ such that $set(L_{a})=f(set(G|a),\vec{y}/\vec{a})$
to $a$ as follows.
If $a$ is the leaf, i.e., the node of outdegree zero,
then $L_{a}$ is a fully collapsed DAG such that
$set(L_{a})=f(\emptyset, \vec{y}/\vec{a})=
h(\emptyset,\vec{y}/\vec{a},\emptyset)$.
Next consider the case when $a$ is not a leaf, and let
$b_{0},\ldots,b_{n}$ be the sons of $a$ in $G$:
$\{b_{0},\ldots,b_{n}\}=\{b\in G: (a,b)\in E_{G}\}$.
Assume that for each son $b_{i}$
a fully collapsed DAG $L_{b_{i}}$ is attached to $b_{i}$ so that
$set(L_{b_{i}})=f(set(G|b_{i}),\vec{y}/\vec{a})$.
Then by Proposition \ref{prp:bisimilarfullcollapsed}.\ref{prp:bisimilarfullcollapsed2}
compute a fully collapsed DAG $C=c(L_{b_{0}},\ldots,L_{b_{n}})$,
and then
let $L_{a}$ be a fully collapsed DAG such that
$set(L_{a})=h(set(G|a),\vec{y}/\vec{a},set(C))$.

Let us estimate roughly the number of computation steps.
The number of number of recursive calls of the function $h$ $h$ is $cT(set(G))+1$.
By Theorem \ref{th:size} we have a polynomial $p_{f}$ such that
\[
cT(L_{a})\leq p_{f}(cT(set(G|a)), cT(set(\vec{H}))) +cT(set(\vec{K}))
.\]
Since all DAG's are balanced and fully collapsed, we have for a polynomial $p_{f}^{\prime}$
\[
|\lceil L_{a}\rceil| \leq p^{\prime}_{f}(|\lceil G|a\rceil|, |\lceil \vec{H}\rceil|, |\lceil\vec{K}\rceil|)
.\]
Hence each computation of $h$ is performed in the number of steps bounded by
a polynomial of $|\lceil G\rceil|$, $|\lceil\vec{H}\rceil|$ and $|\lceil\vec{K}\rceil|$.
Moreover the number of  computations of $C=c(L_{b_{0}},\ldots,L_{b_{n}})$ is $cT(set(G))$,
and  each computation of $C$ is also performed 
polynomially in $|\lceil G\rceil|$, $|\lceil\vec{H}\rceil|$ and $|\lceil\vec{K}\rceil|$.
Hence the number of computation steps for $L$
is bounded by a polynomial of 
$|\lceil G\rceil|$, $|\lceil\vec{H}\rceil|$ and $|\lceil\vec{K}\rceil|$.
\hspace*{\fill} $\Box$

\begin{corollary}\label{th:main}
Suppose a set theoretic function $F(\vec{x})$ is a function on binary finite strings when we restrict to finite strings:
$\forall \vec{s}\subset{}^{<\omega}2\exists t\in{}^{<\omega}2[F(\nu(\vec{s}))=\nu(t)]$.
If
$F(\vec{x}/-)\in{\sf PCSF}$, then
the function $\vec{s}\mapsto \nu^{-1}(F(\nu(\vec{s})))$ is polynomial time computable.
\end{corollary}
\textit{Proof}.\hspace{2mm}
Assume $F\in{\sf PCSF}$, and let $f(\vec{s})=\nu^{-1}(F(\nu(\vec{s})))$.
Then $F$ is a polynomial time function on $\mathbb{HF}$ in the sense of Theorem \ref{lem:hfptime}.
Since the function $s\mapsto \lceil\nu(s)\rceil$ and its inverse $ \lceil\nu(s)\rceil\mapsto s$ are polynomial time computable,
so is $f$.
\hspace*{\fill} $\Box$
\\

\noindent
{\bf Remarks}.
\begin{enumerate}
\item
Let $F$ be a polynomial time computable function for $f\in{\sf PCSF}$ in Theorem \ref{lem:hfptime}.
Then $F$ has to be an `extensional' function on DAG's.
This means that for any balanced and fully collapsed DAG's $\vec{G}$, $\vec{H}$
\[
set(\vec{G})=set(\vec{H}) \,\&\, F(\lceil\vec{G}\rceil)=\lceil K\rceil \,\&\, F(\lceil\vec{H}\rceil)=\lceil L\rceil \Rightarrow
set(K)=set(L)
.\]
\begin{enumerate}
\item
It seems to us that the converse holds.
Namely let $F$ be a polynomial time computable function such that
$F(\lceil\vec{G}\rceil)$ is a code of balanced and fully collapsed DAG for any balanced and fully collapsed DAG's $\vec{G}$, 
and $F$ is extensional in the above sense.
Then the set-theoretic function $f$ on $\mathbb{HF}$ is defined by
$f(\vec{x})=set(H)$ where $\vec{x}=set(\vec{G})$ and $F(\lceil\vec{G}\rceil)=\lceil H\rceil$
for some (any) balanced and fully collapsed DAG's $\vec{G}$ and $H$.

An affirmative answer to the following problem would show a stronger statement than Lemma \ref{lem:ptimePCSF}
since there are polynomial computable functions mapping binary strings $s$ to DAG's (balanced and fully collapsed) 
representing $\nu(s)$, and vice versa.
\\
{\bf Problem}.
Show that the $f$ is a restriction of a function in the class {\sf PCSF} on $\mathbb{HF}$.

\item
Let $c(-/a)$ be a choice function which chooses an element $b\in a$ from non-empty sets $a$.
Let us set $c(-/\emptyset)=\emptyset$.
It is unlikely the case that there is such a $c$ in the class {\sf PCSF}, nor 
$c$ on $\mathbb{HF}$ is (extensionally) polynomial time computable in the sense of Theorem \ref{lem:hfptime}.
Obviously there exists an intensional function $C$ which depends on codes. 
Given DAG's $G$, if $V_{G}\neq\{r_{G}\}$, then let $a_{G}=\max\{a\in V_{G}: a\neq r_{G}\}$.
Then $set(G|a_{G})\in set(G)$, and $\lceil G\rceil \mapsto a_{G}$ is polynomial time computable, and so is
the function $C(\lceil G\rceil)=\lceil G|a_{G}\rceil$.
However $C$ is not extensional.
\end{enumerate}

\item
In \cite{Lago}, U. Dal Lago, S. Martini and M. Zorzi proved that 
ramified recurrence of any free algebra with tiers is polynomial time computable,
as claimed in D. Leivant \cite{Leivant}.
Their proof is based on term graph rewritings, i.e., each term is represented by a term graph (DAG),
and common subterms are sharing.
Our proof of Theorem \ref{lem:hfptime} is akin to their proof in representing data (hereditarily finite sets in our case) as DAG,
but we have to treat a variadic function symbol, i.e., ${\rm Pair}(-/a_{1},\ldots,a_{n})=\{a_{1},\ldots,a_{n}\}$
for $n=0,1,2,3,\ldots$ to represent a hereditarily finite set as a term (graph),
while in \cite{Lago} each function symbol has a fixed arity.
However a term (graph) rewriting approach to {\sf PCSF}-functions is open to us.
\end{enumerate}

\end{document}